\documentclass[12pt]{article}
\usepackage{graphicx}
\usepackage{amsmath,amsthm,amsfonts}
\usepackage{fullpage}

\newtheorem{theorem}{Theorem}
\newtheorem{corollary}[theorem]{Corollary}
\newtheorem{lemma}[theorem]{Lemma}
\newtheorem{proposition}[theorem]{Proposition}
\newtheorem{conjecture}[theorem]{Conjecture}

\title{Words avoiding repetitions in arithmetic progressions}
\author{
Jui-Yi Kao, Narad Rampersad, Jeffrey Shallit \\
David R. Cheriton School of Computer Science \\
University of Waterloo \\
Waterloo, Ontario N2L 3G1 (Canada) \\
{\tt j3kao@student.cs.uwaterloo.ca} \\
{\tt nrampersad@math.uwaterloo.ca} \\
{\tt shallit@graceland.math.uwaterloo.ca} \bigskip \\
Manuel Silva \\
Departamento de Matem\'atica \\
Universidade Nova de Lisboa \\
Quinta da Torre, 2829-516 Caparica (Portugal) \\
{\tt mnas@fct.unl.pt}}

\begin{document}
\date{\today}
\maketitle

\begin{abstract}
Carpi constructed an infinite word over a $4$-letter alphabet
that avoids squares in all subsequences indexed by arithmetic
progressions of odd difference.  We show a connection between
Carpi's construction and the paperfolding words.
We extend Carpi's result by constructing uncountably many words
that avoid squares in arithmetic progressions of odd difference.
We also construct infinite words avoiding overlaps and infinite
words avoiding arbitrarily large squares in arithmetic progressions
of odd difference.  We use these words to construct labelings of the
$2$-dimensional integer lattice such that any line through the
lattice encounters a squarefree (resp.\ overlapfree) sequence of labels.
\end{abstract}

\section{Introduction}

The problem of avoiding repetitions in words was first studied by
Thue \cite{Thu06}, who constructed an infinite word over a ternary
alphabet containing no squares of the form $xx$.
In this paper we generalize this notion by constructing infinite
words containing no squares in any subsequence indexed by an
arithmetic progression of odd difference.  To do so, we make use
of several other generalizations of Thue's problem.

While it is easy to see that any binary word of length at least $4$
must contain a square, Entringer, Jackson, and Schatz \cite{EJS74}
constructed an infinite binary word containing no squares $xx$,
where $|x| \geq 3$.  Prodinger and Urbanek \cite{PU79}
gave an example of an infinite binary word whose only squares
are of lengths $1$, $3$, or $5$.  The particular word studied
by Prodinger and Urbanek is the well-known (ordinary) paperfolding word
\[
0010011000110110\cdots.
\]
Paperfolding words in general have been studied extensively
\cite{All92,AB94b,DMP82}; we will rely in particular on the results of
Allouche and Bousquet-M\'elou \cite{All84,AB94a}.

Taking Thue's problem in another direction, Carpi \cite{Car88},
as a preliminary step in constructing non-repetitive labelings
of the integer lattice, considered the question of the existence
of infinite words that avoid squares in all subsequences
indexed by arithmetic progressions.  Of course, by the classical
theorem of van der Waerden \cite{Wae27}, no such words exist, but Carpi showed
that for any prime $p$, there exists an infinite word over a finite
alphabet that avoids squares in arithmetic progressions of all
differences, except those differences that are a multiple of $p$.  For example,
taking $p = 2$, there exists an infinite word over a $4$-letter alphabet
that contains no squares in any arithmetic progression of odd difference.
As we shall see later, Carpi's construction has a surprising connection
to the paperfolding words.

Another notion of significance in the study of infinite words is that
of \emph{subword complexity}.  The subword complexity function of a
word $w$ is the function $p_w(n)$ that counts the number of distinct
subwords of length $n$ that appear in $w$.  Avgustinovich, Fon-Der-Flaass,
and Frid \cite{AFF03} generalized the concept of subword complexity by
considering the \emph{arithmetical complexity} of a word.  The arithmetical
complexity function of a word $w$ is the function $p^A_w(n)$ that counts the
total number of distinct subwords of length $n$ that appear in all
subsequences of $w$ indexed by arithmetic progressions.  Avgustinovich,
Fon-Der-Flaass, and Frid showed that the words with lowest arithmetical
complexity come from a class of words known as Toeplitz words, of which the
paperfolding words form a special class.  Implicit in their work
is a characterization of the arithmetic subsequences of the
paperfolding words.  We shall rely heavily on this characterization
in our constructions.

\section{Definitions and notation}

Given an infinite word ${\bf w}$ over a finite alphabet $\Sigma$,
we write
\[
{\bf w} = w_0w_1w_2\cdots,
\]
where $w_i \in \Sigma$ for $i \geq 0$.  We sometimes write ${\bf w}[i]$ for
$w_i$.  A \emph{subword} of ${\bf w}$ is a contiguous block of symbols
\[
w_iw_{i+1}\cdots w_{i+j},
\]
for some $i,j \geq 0$.  A \emph{subsequence} of ${\bf w}$ is word of the form
\[
w_{i_0}w_{i_1}\cdots,
\]
where $0 \leq i_0 < i_1 < \cdots$.  An
\emph{arithmetic subsequence of difference $j$} of ${\bf w}$ is a word
of the form
\[
w_i w_{i+j} w_{i+2j} \cdots,
\]
where $i \geq 0$ and $j \geq 1$.  We also define finite subsequences
in the obvious way.

A \emph{square} is a non-empty word $xx$, a \emph{cube} is a non-empty
word $xxx$, and in general, a \emph{$k$-power} is a non-empty word $x^k$.
We define fractional powers in the following way: if $q$ is a positive
rational number, a $q$-power is a non-empty word $x^kx'$, where $x'$ is
a prefix of $x$ and $|x^kx'|/|x| = q$.

If $r$ is a positive real number, we say a word ${\bf w}$
\emph{contains an $r$-power} (resp.\ \emph{contains an $r^+$-power})
if ${\bf w}$ contains a $q$-power as a subword for some $q \geq r$
(resp.\ $q > r$).
A word ${\bf w}$ is \emph{$r$-power-free} (resp.\ \emph{$r^+$-power-free})
or \emph{avoids $r$-powers} (resp.\ \emph{avoids $r^+$-powers})
if ${\bf w}$ contains no $r$-power (resp.\ $r^+$-power).
We use the terms \emph{squarefree}, \emph{overlapfree}, and
\emph{cubefree} for $2$-power-free, $2^+$-power-free, and $3$-power-free,
respectively.

If a word ${\bf w}$ has the property that no arithmetic subsequence of
difference $j$ contains a square (resp.\ cube, $r$-power, $r^+$-power),
we say that ${\bf w}$ \emph{contains no squares (resp.\ cubes, $r$-powers,
$r^+$-powers) in arithmetic progressions of difference $j$}.

For any word $w = w_0w_1\cdots w_n$, we denote by $w^R$ the \emph{reversal}
of $w$, namely the word $w^R = w_nw_{n-1}\cdots w_0$.  For any word $w$
over the binary alphabet $\{0,1\}$, we denote by $\overline{w}$ the
\emph{complement} of $w$, namely the word obtained from $w$ by changing
$0$'s to $1$'s and $1$'s to $0$'s.

\section{Paperfolding words}

A \emph{paperfolding word} ${\bf f} = f_0f_1f_2\cdots$ over the alphabet
$\{0,1\}$ satisfies the following recursive definition:  there exists
$a \in \{0,1\}$ such that
\begin{eqnarray*}
f_{4n} & = & a, \quad n \geq 0 \\
f_{4n+2} & = & \overline{a}, \quad n \geq 0 \\
(f_{2n+1})_{n \geq 0} & \  & \text{is a paperfolding word}.
\end{eqnarray*}
The \emph{ordinary paperfolding word}
\[
0010011000110110\cdots
\]
is the paperfolding word uniquely characterized by $f_{2^m-1} = 0$
for all $m \geq 0$.

One may also define the paperfolding words by means of the
\emph{perturbed symmetry} of Mend\`es France \cite{BM82,Men79}
in the following way.  For $i \geq 0$, let $c_i \in \{0,1\}$
and define the sequence of words
\begin{eqnarray*}
F_0 & = & c_0 \\
F_1 & = & F_0 \,c_1\, \overline{F_0}^R \\
F_2 & = & F_1 \,c_2\, \overline{F_1}^R \\
& \vdots &
\end{eqnarray*}
Then
\[
{\bf f} = \lim_{i\to\infty}F_i
\]
is a paperfolding word.  For example, taking $c_i = 0$ for
all $i \geq 0$, one obtains the sequence
\begin{eqnarray*}
F_0 & = & 0 \\
F_1 & = & 0\,0\,1 \\
F_2 & = & 001\,0\,011 \\
& \vdots &
\end{eqnarray*}
which converges, in the limit, to the ordinary paperfolding word.

The following properties of paperfolding words were
proved by Allouche and Bousquet-M\'elou \cite{All84,AB94a}
(the particular case of the ordinary paperfolding word
was studied by Prodinger and Urbanek \cite{PU79}).

\begin{theorem}[Allouche and Bousquet-M\'elou]
\label{all1}
For any paperfolding word ${\bf f}$, if $xx$ is a non-empty
subword of ${\bf f}$, then $|x| \in \{1,3,5\}$.
\end{theorem}

\begin{corollary}[Allouche and Bousquet-M\'elou]
\label{all2}
For any paperfolding word ${\bf f}$, ${\bf f}$ contains no
fourth powers and no cubes except $000$ and $111$.  In particular,
${\bf f}$ contains no $3^+$-power.
\end{corollary}

Unfortunately, the proof of Theorem~\ref{all1} given in
\cite{AB94a} contains an error.  For completeness we therefore
provide a proof below.  We first prove the following corrected version of
\cite[Proposition~5.1]{AB94a}.

\begin{proposition}
\label{fixed}
If a paperfolding word ${\bf f}$ contains a subword $wcw$, where $w$
is a non-empty word and $c$ is a single letter, then
either $|w| \in \{2,4\}$ or $|w| = 2^k - 1$ for some $k \geq 1$.
\end{proposition}

We will need the following result due to Allouche \cite{All92}.

\begin{lemma}[Allouche]
\label{evenodd}
Let $u$ and $v$ be subwords of a paperfolding word ${\bf f}$,
with $|u| = |v| \geq 7$.  If $u$ and $v$ occur at positions
of different parity in ${\bf f}$, then $u \neq v$.
			     \end{lemma}

\begin{proof}[Proof of Proposition~\ref{fixed}]
Suppose to the contrary that
\[
wcw = f_i f_{i+1} \cdots f_{i+t} f_{i+t+1} \cdots f_{i+2t}
\]
is a subword of ${\bf f}$, where $|w| = t$, $t \notin \{2,4\}$,
$t \neq 2^k - 1$ for all $k \geq 1$.  Suppose further that ${\bf f}$
is chosen so as to minimize $t$.  We consider four cases.

Case 1: $t = 6$.  Because the letters in successive even positions of ${\bf f}$
alternate between $0$ and $1$, any subword of ${\bf f}$ of length $13$ starting
at an even position must be of the form
\begin{eqnarray*}
0*1*0*1*0*1*0 & \text{or} & 1*0*1*0*1*0*1,
\end{eqnarray*}
where the $*$ denotes an arbitrary symbol from $\{0,1\}$.  Consequently,
if such a subword is of the form $wcw$, it must be one of the words
\begin{eqnarray*}
0011001001100 & \text{or} & 1100110110011.
\end{eqnarray*}
Similarly, if $wcw$ begins at an odd position, it must be one of the words
\begin{eqnarray*}
011001c011001 & \text{or} & 100110c100110.
\end{eqnarray*}
Taking the odd indexed positions of $wcw$, we see that
if $i$ is even, then either $010010$ or $101101$ is a subword of a paperfolding
word, which is impossible, since neither word obeys the required alternation of
$0$'s and $1$'s in even indexed positions.  Similarly, if $i$ is odd,
then either $010c101$ or $101c010$ is a subword of a paperfolding word,
which again is impossible for any choice of $c$.

Case 2: $t$ even, $t \geq 8$.  Then $w$ occurs at positions of two different
parities in ${\bf f}$, contradicting Lemma~\ref{evenodd}.

Case 3: $t \equiv 1 \pmod 4$, $t \geq 5$.  Let $\ell \in \{i,i+1\}$
such that $\ell$ is even.  Then $f_{\ell} \neq f_{\ell+t+1}$,
since $\ell$ and $\ell+t+1$ are even but $\ell \not\equiv \ell+t+1 \pmod 4$.

Case 4: $t \equiv 3 \pmod 4$, $t \geq 11$.  Let $t = 4m+3$, where $m \geq 2$
and $m+1$ is not a power of $2$.  Let $\ell \in \{i,i+1\}$ such
that $\ell$ is odd.  Then
\[
w'c'w' = f_{\ell} f_{\ell+2} \cdots f_{\ell+t-1} f_{\ell+t+1} \cdots
f_{\ell+2t-2}
\]
is a subword of a paperfolding word, where $|w'| = t' = (t-1)/2 = 2m+1$.
By the argument of Case~3, $t' \not\equiv 1 \pmod 4$.  Let us write
$t' = 4m'+3$, where $m' = (m-1)/2$.  Since $m+1$ is not a power of $2$,
$m'+1$ is not a power of $2$.  Thus $11 \leq t' < t$, contradicting the
minimality of $t$.
\end{proof}

The following result is not needed for the proof of Theorem~\ref{all1}
but will be useful in the next section.

\begin{proposition}
\label{wcw}
Let ${\bf f}$ be a paperfolding word.  For all $k \geq 1$, ${\bf f}$
contains a subword $wcw$, where $w$ is a non-empty word, $c$ is a single
letter, and $|w| = 2^k-1$.
\end{proposition}

\begin{proof}
By the perturbed symmetry construction, ${\bf f}$ begins with a prefix
$z c_0 \overline{z}^R$, where $|z| = 2^{k-1}-1$ and $c_0 \in \{0,1\}$.
Applying the perturbed symmetry map twice to $z c_0 \overline{z}^R$,
we see that ${\bf f}$ begins with a prefix
\[
z \;c_0\; \overline{z}^R \;c_1\; z \;\overline{c_0}\; \overline{z}^R \;c_2\;
z \;c_0\; \overline{z}^R \;\overline{c_1}\; z \;\overline{c_0}\;
\overline{z}^R,
\]
where $c_1,c_2 \in \{0,1\}$.  If $c_1 = c_2$, then
\[
wcw = \overline{z}^R \;c_1\; z \;\overline{c_0}\; \overline{z}^R \;c_2\; z
\]
is the desired subword.
If $c_1 \neq c_2$, then
\[
wcw = \overline{z}^R \;c_2\; z \;c_0\; \overline{z}^R \;\overline{c_1}\; z
\]
is the desired subword.
\end{proof}

We will also need the following lemma.

\begin{lemma}
\label{powersof2}
For all $k \geq 1$, no paperfolding word ${\bf f}$ contains a subword $xx$
with $|x| = 2^k$.
\end{lemma}

\begin{proof}
The proof is by induction on $k$.  If $k=1$, then let
$f_i f_{i+1} f_{i+2} f_{i+3}$ be a subword of ${\bf f}$.
If $i$ is even (resp.\ odd), then $f_i \neq f_{i+2}$ (resp.\
$f_{i+1} \neq f_{i+3}$).

Now suppose
\[
xx = f_i f_{i+1} \cdots f_{i+2^{k+2}-1}
\]
is a subword of ${\bf f}$.  Let $\ell \in \{i,i+1\}$ such that $\ell$ is odd.
Then
\[
x'x' = f_{\ell} f_{\ell+2} \cdots f_{\ell+2^{k+2}-2}
\]
is a subword of a paperfolding word with $|x'| = 2^k$.  The result follows
by induction.
\end{proof}

We are now ready to prove Theorem~\ref{all1}.

\begin{proof}[Proof of Theorem~\ref{all1}]
If ${\bf f}$ contains a square $xx$, then writing $x = wc$, where $c$ is a
single letter, we see that ${\bf f}$ contains the subword $wcw$.  By
Proposition~\ref{fixed}, either $|x| \in \{1,3,5\}$, or $|x| = 2^k$
for some $k \geq 1$.  But we have seen in Lemma~\ref{powersof2}
that the latter is impossible.
\end{proof}

We end this section with the following interesting fact regarding the
ordinary paperfolding word.

\begin{proposition}
Let ${\bf f}$ be the ordinary paperfolding word over $\{0,1\}$.
Then $0{\bf f}$ is the lexicographically least word in the orbit
closure of any paperfolding word.
\end{proposition}

\begin{proof}
Taking the subsequence of ${\bf f}$ indexed by the odd positions
yields the word ${\bf f}$ again, so taking the subsequence of
$0{\bf f}$ indexed by the even positions yields the word $0{\bf f}$.

Let ${\bf w} = w_0w_1w_2\cdots$ be the lexicographically least word
in the orbit closure of any paperfolding word.  Let us assume that
${\bf w}$ begins with $0001$, since it cannot begin with anything
lexicographically smaller.  Since $w_0 = w_2$, the following is
forced:  $w_1w_3w_5w_7\cdots = 0101\cdots$.

We will prove by induction on $n$ that the prefixes of ${\bf w}$
of length $2n$ are the prefixes of $0{\bf f}$.  We have already
established the base case, so let us suppose $n \geq 2$ and
$w_0w_1w_2\cdots w_{2n-1} = 0f_0f_1f_2\cdots f_{2n-2}$.  Since
$w_1w_3w_5w_7\cdots = 0101\cdots$, we see that $w_{2n+1} = f_{2n}$.
Note that $w_0w_2w_4\cdots w_{2n} = 0f_1f_3f_5\cdots f_{2n-1}$
is a prefix of a word in the orbit closure of a paperfolding word.
By our inductive assumption, $w_0w_1w_2\cdots w_{n-1}w_n$ is
the lexicographically least such prefix.  Choosing
$w_{2n} = w_n = f_{n-1} = f_{2n-1}$ thus ensures that
$w_0w_1w_2\cdots w_{2n}w_{2n+1}$ is lexicographically minimal.
We have thus established that ${\bf w}$ and $0{\bf f}$ agree on
the first $2(n+1)$ positions, as required.
\end{proof}

\section{Avoiding repetitions in arithmetic progressions}

In this section we construct infinite words avoiding squares (resp.\
overlaps) in all arithmetic progressions of odd difference.

The following result is implicit in the work of Avgustinovich,
Fon-Der-Flaass, and Frid (see the proof of \cite[Theorem~3]{AFF03}
as well as \cite[Example~2]{AFF03}).

\begin{theorem}[Avgustinovich, Fon-Der-Flaass, and Frid]
\label{frid}
If $w$ is a finite arithmetic subsequence of odd difference of
a paperfolding word, then $w$ is a subword of a paperfolding word.
\end{theorem}

\begin{corollary}
\label{threeplus}
There exists an infinite word over a binary alphabet that contains
no $3^+$-powers in arithmetic progressions of odd difference.
\end{corollary}

\begin{proof}
It follows from Corollary~\ref{all2} and Theorem~\ref{frid} that all
paperfolding words have this property.
\end{proof}

We note further that the $3^+$ of the preceding corollary may not
be replaced by $3$.  The usual backtracking search suffices to verify that
all sufficiently long binary words contain a cube in an arithmetic
progression of odd difference.  The longest binary words that do
not contain a cube in an arithmetic progression of odd difference
are the following words of length $13$:
\begin{eqnarray*}
0010011001100 & \quad & 0101100110011 \\
1010011001100 & \quad & 1101100110011.
\end{eqnarray*}

The problem of avoiding repetitions in arithmetic progressions
seems to have first been studied by Carpi \cite{Car88} and
subsequently by Currie and Simpson \cite{CS02}.  Downarowicz
\cite{Dow99} studied a related problem.

\begin{theorem}[Carpi]
\label{carpi1}
There exists an infinite word over a $4$-letter alphabet that contains
no squares in arithmetic progressions of odd difference.
\end{theorem}

The word ${\bf c}$ constructed by Carpi satisfying the conditions of this
theorem is over the alphabet $\{1,3,5,7\}$ and is generated by
iterating the morphism $1 \to 53$, $3 \to 73$, $5 \to 51$, $7 \to 71$,
starting with the symbol $5$.  It can also be derived from a paperfolding
sequence, as we shall see below.  The alphabet size of $4$ in
Theorem~\ref{carpi1} is optimal, since the longest words over
the alphabet $\{0,1,2\}$ that avoid squares in all odd difference arithmetic
progressions are the words
\[
010212021 \quad 012010201
\]
of length $9$, along with the words obtained from these by permuting
the alphabet symbols in all possible ways.

Let ${\bf f} = f_0f_1f_2\cdots$ be any paperfolding word over
$\{1,4\}$.  Define ${\bf v} = v_0v_1v_2\cdots$ by
\begin{eqnarray*}
v_{4n} & = & 2 \\
v_{4n+2} & = & 3 \\
v_{2n+1} & = & f_{2n+1},
\end{eqnarray*}
for all $n \geq 0$.
In other words, we have recoded the periodic subsequence formed
by taking the even positions of ${\bf f}$ by mapping $1 \to 2$
and $4 \to 3$ (or vice-versa).  For example, if
\[
{\bf f} = 1141144111441441\cdots
\]
is the ordinary paperfolding word over $\{1,4\}$, then
\[
{\bf v} = 2131243121342431\cdots.
\]

\begin{theorem}
\label{squares}
Let ${\bf v}$ be any word obtained from a paperfolding word ${\bf f}$
by the construction described above.  Then the word ${\bf v}$
contains no squares in arithmetic progressions of odd
difference but does not avoid $r$-powers for any real $r<2$.
\end{theorem}

\begin{proof}
By the construction of ${\bf v}$, any arithmetic
subsequence
\[
w = v_{i_0}v_{i_1}\cdots v_{i_k}\]
of odd difference of ${\bf v}$ can be obtained from the
corresponding subsequence
\[
x = f_{i_0}f_{i_1}\cdots f_{i_k}
\]
of ${\bf f}$ by recoding the symbols in
either the even positions of $x$ or the odd positions of $x$ by
mapping  $1 \to 2$ and $4 \to 3$ (or vice-versa).
Note that this recoding cannot create any new squares.  Now suppose that
${\bf v}$ contains a square $ww$ in an arithmetic progression of
odd difference.  Let $xx$ be the corresponding subsequence of ${\bf f}$.
By Theorems~\ref{all1} and~\ref{frid}, $|x| \in \{1,3,5\}$ and hence
$|w| \in \{1,3,5\}$.  Clearly, $|w| = 1$ is impossible.
If $|w| = 3$, then $ww$ has one of the forms
$({}*2*{})(\,3*2\,)$, $({}*3*{})(\,2*3\,)$, $(\,2*3\,)({}*2*{})$, or
$(\,3*2\,)({}*3*{})$, where the $*$ denotes an arbitrary symbol
from $\{1,4\}$.  Clearly, none of these can be squares.
A similar argument applies for $|w| = 5$.

That ${\bf v}$ does not avoid $r$-powers for any $r<2$ follows easily
from Proposition~\ref{wcw}.
\end{proof}

The word ${\bf c}$ constructed by Carpi, after relabeling the alphabet symbols
by the map $1 \to 2$, $3 \to 3$, $5 \to 1$, $7 \to 4$, is the word $1{\bf v}$,
where ${\bf v}$ is constructed from the ordinary paperfolding word as
described above.  Note that since there are uncountably many
paperfolding words ${\bf f}$, there are uncountably many words ${\bf v}$
over a $4$-letter alphabet that contain no squares in arithmetic progressions
of odd difference.  We offer the following conjectures regarding
such words.

\begin{conjecture}
\label{lt2}
For all real numbers $r < 2$, $r$-powers are not avoidable in
arithmetic progressions of odd difference over a 4-letter alphabet.
\end{conjecture}

A backtracking search confirms that Conjecture~\ref{lt2} holds
for all $r \leq 7/4$.

\begin{conjecture}
Any infinite word over a 4-letter alphabet that avoids squares
in arithmetic progressions of odd difference is in the orbit closure
of one of the words ${\bf v}$ constructed above.
\end{conjecture}

Next we consider words over a ternary alphabet.

\begin{theorem}
\label{overlaps}
There exists an infinite word over a ternary alphabet that contains
no $2^+$-powers (overlaps) and no squares $xx$, $|x| \geq 2$,
in arithmetic progressions of odd difference.
\end{theorem}

\begin{proof}
Let ${\bf v} = v_0v_1v_2\cdots$ be any word obtained from a
paperfolding word by the construction described above.
Let $h$ be the morphism that sends $1 \to 00$, $2 \to  11$, $3 \to 12$,
$4 \to 02$.  Then ${\bf w} = w_0w_1w_2\cdots = h({\bf v})$ has the desired
properties.

Suppose to the contrary that there exists $i \geq 0$, $j$ odd, and $t \geq 2$
such that for $s \in \{0,\ldots,t-1\}$, $w_{i+sj} = w_{i+(s+t)j}$.
Note that there exists $a \in \{0,1\}$ such that for $\ell \equiv 0 \pmod 4$,
$w_\ell = a$ and for $\ell \equiv 2 \pmod 4$, $w_\ell = \overline{a}$.
We consider four cases.

Case 1: $t = 3$.  Because the letters in successive even positions of ${\bf w}$
alternate between $0$ and $1$, $w_i w_{i+j} \cdots w_{i+5j}$
is one of the words $001001$, $011011$, $100100$, or $110110$.
Thus there exists $s \in \{0,1,2\}$ such that
$w_{i+sj} \neq w_{i+(s+4)j}$. Now consider the morphism $h$.
The symbol $0$ only occurs in the images of $1$ and $4$, and
the symbol $1$ only occurs in the images of $2$ and $3$.
Let $i' = \lfloor (i+sj)/2 \rfloor$.
Since $w_{i+sj} \neq w_{i+(s+4)j}$, we have that either
$v_{i'} \in \{1,4\}$ and $v_{i'+2j} \in \{2,3\}$, or vice versa.
Either case is impossible, since the symbols $1$ and $4$
only occur in positions of odd parity in ${\bf v}$,
and the symbols $2$ and $3$ only occur in positions of
even parity in ${\bf v}$, but $i'$ and $i'+2j$ both
have the same parity.

Case 2: $t$ odd, $t \geq 5$.  Since $j$ is odd,
$\{i \pmod 8, i+j \pmod 8, \ldots, i + (2t-1)j \pmod 8\}$ is a complete
set of residues $\pmod 8$.  Since ${\bf v}$ contains a $3$ in every position
congruent to $2 \pmod 4$, ${\bf w}$ contains a $2$ in every
position congruent to $5 \pmod 8$.  Thus there exists
$s \in \{0,\ldots,2t-1\}$ such that $w_{i+sj} = 2$.  If $s < t$,
then since $t$ is odd, $s \not\equiv s+t \pmod 2$, and consequently,
$i+sj \not\equiv i+(s+t)j \pmod 2$.  But ${\bf w}$
only contains $2$'s in positions of even parity,
so $w_{i+sj} \neq w_{i+(s+t)j}$, contrary to our assumption.
Similarly, if $s \geq t$, we have $w_{i+(s-t)j} \neq w_{i+sj}$.

Case 3: $t \equiv 2 \pmod 4$.  Then either $w_i \neq w_{i+tj}$ or
$w_{i+j} \neq w_{i+(t+1)j}$, accordingly as $i$ is even or odd,
contrary to our assumption.

Case 4: $t \equiv 0 \pmod 4$.  Let $k \in \{i,i+j\}$ such that $k$ is odd.
Let $k' = \lfloor k/2 \rfloor$.  It follows from the definition of $h$ that
for $s \in \{0,\ldots,t-1\}$, $v_{k'+sj}$ is uniquely
determined by the value of $w_{k+2sj}$ and the congruence class
of $k+2sj \pmod 4$:
\begin{itemize}
\item if $w_{k+2sj} = 0$, then $v_{k'+sj} = 1$;
\item if $w_{k+2sj} = 1$, then $v_{k'+sj} = 2$; and
\item if $w_{k+2sj} = 2$, then $v_{k'+sj}$ is either $3$ or $4$,
accordingly as $k+2sj \equiv 1$ or $3 \pmod 4$.
\end{itemize}
From this observation, combined with the fact that
$k+2sj \equiv k+(2s+t)j \pmod 4$, we see that
since $w_k w_{k+2j} \cdots w_{k+2(t-1)j}$ is a square,
$v_{k'} v_{k'+j} \cdots v_{k'+(t-1)j}$ is also
a square in an arithmetic progression of odd difference $j$ in ${\bf v}$,
a contradiction.

These four cases cover all possibilities.  It remains to consider
the existence of the cubes $000$, $111$, and $222$.  Suppose there
exists $w_i w_{i+j} w_{i+2j} \in \{000,111,222\}$ for some $i \geq 0$
and $j$ odd.  Since ${\bf w}$ only contains $2$'s in positions
of even parity, we may suppose $w_i w_{i+j} w_{i+2j} \in \{000,111\}$.
If $i$ is even, then $i+2j$ is even and $i \not\equiv i+2j \pmod 4$, so
$w_i \neq w_{i+2j}$.  If $i$ is odd, then by the same reasoning
as in Case 4 above, $v_{\lfloor i/2 \rfloor} v_{\lfloor i/2 \rfloor + j}$
is a square in an arithmetic progression of odd difference in ${\bf v}$,
a contradiction.
\end{proof}

The alphabet size of $3$ in Theorem~\ref{overlaps} is optimal, since the
longest words over the alphabet $\{0,1\}$ that avoid overlaps in all odd
difference arithmetic progressions are the words
\[
0010011001 \quad 0101100110 \quad 0110100101
\]
of length $10$, along with their complements.

\section{Avoiding arbitrarily large squares}

In this section we improve upon the result of Entringer, Jackson, and Schatz
\cite{EJS74} noted in the introduction.

\begin{theorem}
\label{largesqs}
There exists an infinite word over a binary alphabet that contains
no squares $xx$ with $|x| \geq 3$ in any arithmetic progression of odd
difference.
\end{theorem}

\begin{proof}
Let ${\bf v}$ be any word obtained from a
paperfolding word by the construction described in the previous section.
Let $h$ be the morphism that sends
\begin{eqnarray*}
1 & \to & 0110 \\
2 & \to & 0101 \\
3 & \to & 0001 \\
4 & \to & 0111.
\end{eqnarray*}
We will show that $h({\bf v})$ has the desired properties.
We first proceed to prove two lemmas about $h({\bf v})$.

\begin{lemma}
\label{largesqs-partialcharac}
Every finite subword \(\alpha\) of an arithmetic subsequence of odd difference
of \(h({\bf v})\) is also a subword of \({\bf W} = \prod_{i\geq 0} W_i \),
where \({\bf W}\) satisfies one of the following conditions:

\begin{itemize}
\item[(a)] \(W_i \in \{0011, 0111\}\) when \(i\) is odd and
\(W_i \in \{0100, 0101\}\) when \(i\) is even.

\item[(b)] \(W_i \in \{0110, 0111\}\) when \(i\) is odd and
\(W_i \in \{0101, 0001\}\) when \(i\) is even.
\end{itemize}
\end{lemma}

\begin{proof}
Any finite subsequence $\alpha$ is a subword of an infinite
subsequence ${\bf W} = (h({\bf v})[q+id])_{i \geq 0}$, where
$q \in \{0,1,2,3\}$ and $d$ is odd.  We have four cases for $d$, namely,
$d \equiv 1,3,5$ or $7 \pmod 8$, respectively.

Suppose $d \equiv 1 \pmod 8$.  Let us also take $q = 0$.
It will be clear from what follows that we may do this with no loss
of generality.  The sequence
\[
(id \bmod 4)_{i \geq 0} = 0,1,2,3,0,1,2,3,\ldots
\]
is periodic with period $4$, and the sequence
\[
\left(\left\lfloor id/4 \right\rfloor \bmod 2\right)_{i \geq 0} =
0,0,0,0,1,1,1,1,0,0,0,0,1,1,1,1,\ldots
\]
is periodic with period $8$.  Note that for $\lfloor id/4 \rfloor \equiv
0 \pmod 2$, ${\bf v}[\,\lfloor id/4 \rfloor\,] \in \{2,3\}$, and for
$\lfloor id/4 \rfloor \equiv 1 \pmod 2$, ${\bf v}[\,\lfloor id/4 \rfloor\,]
\in \{1,4\}$.  Since $h(2)$ and $h(3)$ are equal in all but the second
position, we see that for $\lfloor id/4 \rfloor \equiv 0 \pmod 2$ and
$i \equiv 0 \pmod 4$, we have
\[
(h({\bf v})[(i+j)d])_{j=0,1,2,3} \in \{0101,0001\}.
\]
Similarly, since $h(1)$ and $h(4)$ are equal in all but the last position,
we see that for $\lfloor id/4 \rfloor \equiv 1 \pmod 2$ and
$i \equiv 0 \pmod 4$, we have
\[
(h({\bf v})[(i+j)d])_{j=0,1,2,3} \in \{0110,0111\}.
\]
Thus ${\bf W}$ satisfies condition~(b), as required.
The analysis for $d \equiv 7 \pmod 8$ is similar and results
in ${\bf W}$ satisfying condition~(a).

Now suppose $d \equiv 5 \pmod 8$.  Again we take $q = 0$.
The argument is similar to that for $d \equiv 1 \pmod 8$, except
we consider the sequences
\[
(id \bmod 4)_{i \geq 0} = 0,1,2,3,0,1,2,3,\ldots
\]
and
\[
\left(\left\lfloor id/4 \right\rfloor \bmod 2\right)_{i \geq 0} =
0,1,0,1,1,0,1,0,0,1,0,1,1,0,1,0,\ldots,
\]
where the latter is again periodic with period $8$.  In this case
we deduce that ${\bf W}$ satisfies condition~(a).
The analysis for $d \equiv 3 \pmod 8$ is similar and results
in ${\bf W}$ satisfying condition~(b).
\end{proof}

\begin{lemma}
\label{largesqs-easy_cases}
The word \(h({\bf v})\) contains no squares $xx$ with $|x| = 4$ or $|x| \geq 3$
and $|x| \not\equiv 0 \pmod 4$.
\end{lemma}

\begin{proof} 
Suppose to the contrary that $h({\bf v})$ contains such a square $xx$.
Let $xx$ be a subword of $\prod_{i \geq 0} W_i$, as in
Lemma~\ref{largesqs-partialcharac}.
We consider five cases.  In Cases~1--3, let $xx$ be a subword of
$W_{q} \cdots W_{q+2k}$ for some $q$ and some minimal $k$.  Let us also
write
\[
W_{q} \cdots W_{q+2k} = A_0 A_1 \cdots A_k B_1 \cdots B_k,
\]
where for $i = 0,\ldots,k$, $A_i = W_{q+i}$ and for $i = 1,\ldots,k$,
$B_i = W_{q+k+i}$.  We also define $B_0 = A_k$.

Case 1: $|x| \equiv 1 \pmod 4$ and $|x| \geq 9$. The situation is depicted in Figure~\ref{fig:1mod4}. It is clear from the figure that $A_1[0] = B_1[1]$ and $A_2[0] = B_2[1]$. But from Lemma~\ref{largesqs-partialcharac}, $A_1[0] = 0 = A_2[0]$. Checking the two conditions given in Lemma~\ref{largesqs-partialcharac} shows that $B_1[1] = B_2[1] = 0$ is a contradiction.

\begin{figure}[ht]
\centering
\includegraphics{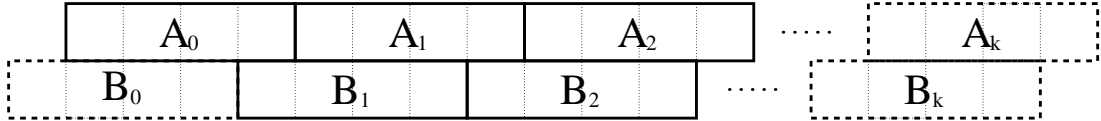}
\caption{$|x| \equiv 1 \pmod 4$ and $|x| \geq 9$}
\label{fig:1mod4}
\end{figure}

Case 2: $|x| \equiv 2 \pmod 4$ and $|x| \geq 9$. The situation is depicted in Figure~\ref{fig:2mod4}. It is clear from the figure that $A_1[0] = B_1[2]$ and $A_2[0] = B_2[2]$. But from Lemma~\ref{largesqs-partialcharac}, $A_1[0] = 0 = A_2[0]$. Checking the two conditions given in Lemma~\ref{largesqs-partialcharac} shows that $B_1[2] = B_2[2] = 0$ is a contradiction.

\begin{figure}[ht]
\centering
\includegraphics{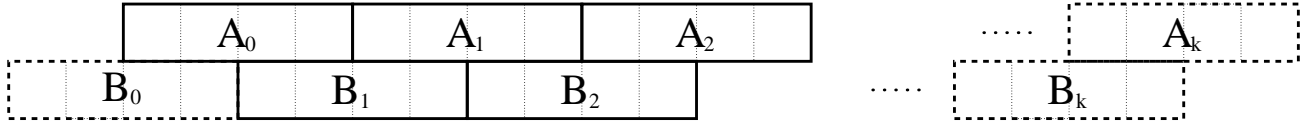}
\caption{$|x| \equiv 2 \pmod 4$ and $|x| \geq 9$}
\label{fig:2mod4}
\end{figure}

Case 3: $|x| \equiv 3 \pmod 4$ and $|x| \geq 9$. The situation is depicted in Figure~\ref{fig:3mod4}. It is clear from the figure that $A_1[0] = B_1[3]$ and $A_2[0] = B_2[3]$. But from Lemma~\ref{largesqs-partialcharac}, $A_1[0] = 0 = A_2[0]$. Checking the two conditions given in Lemma~\ref{largesqs-partialcharac} shows that $B_1[3] = B_2[3] = 0$ is a contradiction.

\begin{figure}[ht]
\centering
\includegraphics{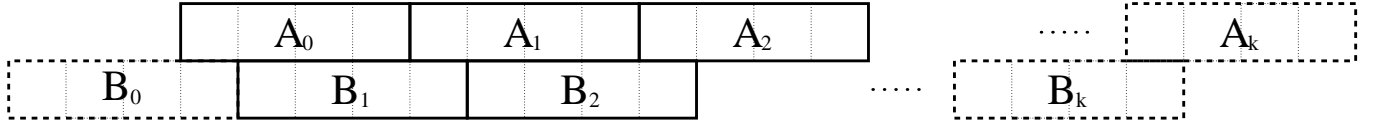}
\caption{$|x| \equiv 3 \pmod 4$ and $|x| \geq 9$}
\label{fig:3mod4}
\end{figure}

Case 4: $|x| = 3,4,5$ or $6$. Let $xx$ be a subword of $A_0 A_1 A_2 A_3$ where for some $p$ and for each $i = 0,1,2,3$, $A_i = W_{p+i}$. By Lemma~\ref{largesqs-partialcharac}, there are at most $64$ possibilities for $A_0 A_1 A_2 A_3$. It is easy to check with the aid of a computer that none of these words contain squares of length greater than $3$.

Case 5: $|x| = 7$.  Let $xx$ be a subword of $A_0 A_1 A_2 A_3 A_4$ where for some $p$ and for each $i = 0,1,2,3,4$, $A_i = W_{p+i}$. For some $q \in \{0,1,2,3\}$, $(xx)[i] = (A_0 A_1 A_2 A_3 A_4)[q+i]$ for all $i \in \{0,\ldots,2|x|-1\}$.
If $q \in \{0,1,2\}$, then $A_4$ is irrelevant. Case~4 above shows that
no such square occurs. Otherwise, $q = 3$. We then have
\[
x = A_0[3] \, A_1\, A_2[0]\, A_2[1] = A_2[2] \, A_2[3]\, A_3\, A_4[0].
\]
In particular, $A_2[1] = A_4[0]$ and $A_2[3] = A_1[0]$.  Since $W_i[0] = 0$
for all $i \geq 0$, we have $A_2[0] = A_2[1] = A_2[3] = 0$.
There is no such $W_i = A_2$ by Lemma~\ref{largesqs-partialcharac}.
\end{proof}

To complete the proof of Theorem~\ref{largesqs}, it remains to consider the
case where \(|x| \equiv 0 \pmod 4\), \(|x| \geq 8\).  Suppose that for such
an $x$, \( xx\) occurs as an arithmetic subsequence of odd difference in
\(h({\bf v})\).

Let \(y \in \{1,2,3,4\}^*\) and \(z = h(y) \in \{0,1\}^*\) such that \(y\) is a minimal subword of \({\bf v}\) such that \(xx\) occurs over an odd-difference arithmetic progression over \(z = h(y)\).
That is, for some fixed \(q \in \{0,1,2,3\}\) and \(d\) a positive odd integer, \( xx = (z[q+id])_{i = 0,..,2|x|-1} \).
We will derive a contradiction by showing that $y$ contains a square in an odd-difference arithmetic progression.

Let $l \in \{0,1,2,3\}$ such that $y[l] = 3$.  Since $d$ is odd,
one easily verifies that there exists $i_0$, $0 \leq i_0 \leq 15$,
satisfying $q+i_0 d \equiv 1 \pmod 4$ and
$\lfloor (q+i_0 d)/4 \rfloor \equiv l \pmod 4$, so that
$y[ \, \lfloor (q+i_0 d)/4 \rfloor \, ] = 3$.
Fix such an \(i_0\). If \(i_0 \in \{0,\ldots,|x|-1\}\), then
\[
z[q+i_0 d] = z[q+(|x|+i_0)d] = 0.
\]
Since \(|x| \equiv 0 \pmod 4\), we have
\[
q+i_0 d \equiv q+(|x|+i_0)d \equiv 1 \pmod 4,
\]
so
\[
h(y[ \, \lfloor (q+(|x|+i_0) d)/4 \rfloor \, ])[1] = 0.
\]
A quick check of the possible images of \(h\) shows that
\(y[ \, \lfloor (q+(|x|+i_0) d)/4 \rfloor \, ] = 3\).

Similarly, if \(i_0 \in \{|x|,\ldots,2|x|-1\}\), then
\(i_0 - |x| \in \{0,\ldots,|x|-1\}\) satisfies the same requirements.
Without loss of generality, we may assume \(i_0 \in \{0,\ldots,|x|-1\}\).

\begin{figure}[htb]
\centering
\includegraphics{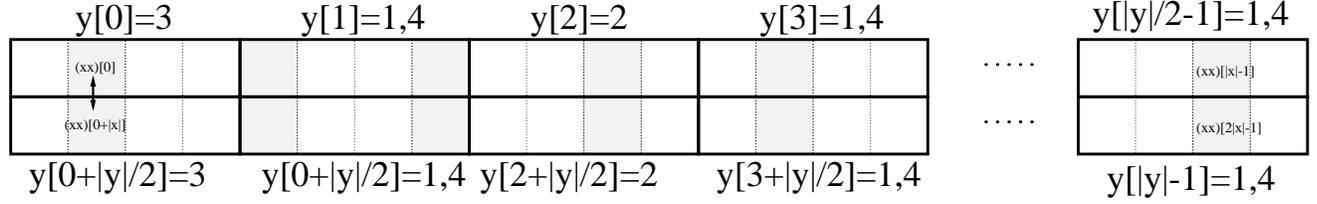}
\caption{An example illustrating the characterization of $y$}
\label{fig:0mod4-1}
\end{figure}

Let $b_1,b_2 \in \{0,1\}$ be such that $y[b_1] \in \{2,3\}$ and
$y[b_1 + 2b_2] = 3$.
Then we can characterize \(y\) as follows (Figure~\ref{fig:0mod4-1}):
for $j \in \{0,\ldots,|y|-1\}$

\begin{itemize}
\item[(a)] If \(j \equiv b_1 +2 b_2 \pmod 4\), then \(y[j] = y[j+|y|/2] = 3\).

\item[(b)] If \(j \equiv b_1 +2 (b_2 + 1) \pmod 4\), then \(y[j] = y[j+|y|/2] = 2\).

\item[(c)] If \(j \not\equiv b_1 \pmod 2\), then \(y[j],y[j+|y|/2] \in \{1,4\}\).
\end{itemize}

Consider the simultaneous congruences
\begin{equation*}
\begin{cases}
s \equiv q \pmod d; \\
s \equiv 3 \pmod 4.
\end{cases}
\end{equation*}

The solution is of the form \(s = s_0 + m \cdot 4d\) for all \(m\), where
\(s_0\) is the least solution in the range \(\{q,\ldots,q+(|x|-1)d\}\). 
Let \(m_0\) be such that \(s_0 + m_0 \cdot 4d\) is the greatest solution
in the range \(\{q,\ldots,q+(|x|-1)d\}\).
Consider each \(s = s_0 + m \cdot 4d\) in the range
\(\{q,\ldots,q+(|x|-1)d\}\).
If \(j = \lfloor s/4 \rfloor \equiv b_1 (\)mod \(2)\), then
\(y[j] = y[j+|y|/2]\) by (a) and (b) above.
If \(j = \lfloor s/4 \rfloor \not\equiv b_1 (\)mod \(2)\), then by (c)
\(y[j]\) and \(y[j+|y|/2] \in \{1,4\}\) and \(h(y[j])[3] = h(y[j+|y|/2])[3]\).
Since \(h(1)[3] \neq h(4)[3]\), we have \(y[j] = y[j+|y|/2]\).

Let $c = \lfloor s_0/4 \rfloor$.  Then
\begin{multline*}
(y[ \, \lfloor (s_0 + m \cdot 4d)/4 \rfloor \, ])_{m = 0,\ldots,m_0} = \\
y[c] \, y[c+d] \, y[c+2d] \, \cdots \,  y[c+m_0 d] \,
y[c+|y|/2] \, y[c+d+|y|/2] \\ y[c+2d+|y|/2] \, \cdots \,  y[c+m_0 d+|y|/2]
\end{multline*}
is a square in an odd-difference arithmetic progression over \({\bf v}\),
contradicting Theorem~\ref{squares}.
\end{proof}

\section{Avoiding repetitions in higher dimensions}

An infinite word ${\bf w}$ over a finite alphabet $A$ is a map
from $\mathbb{N}$ to $A$, where we write $w_n$ for ${\bf w}(n)$.
Now consider a map ${\bf w}$ from $\mathbb{N}^2$ to $A$,
where we write $w_{m,n}$ for ${\bf w}(m,n)$.  We call such a ${\bf w}$
a \emph{2-dimensional word}.
A word ${\bf x}$ is a \emph{line} of ${\bf w}$ if there exists
$i_1,i_2$, $j_1,j_2$ such that $\gcd(j_1,j_2) = 1$, and for $t \geq 0$,
\[
x_t = w_{i_1+j_1t,i_2+j_2t}.
\]

Carpi \cite{Car88} proved the following surprising result.

\begin{theorem}[Carpi]
\label{carpi2}
There exists a 2-dimensional word ${\bf w}$ over a 16-letter alphabet,
such that every line of ${\bf w}$ is squarefree.
\end{theorem}

\begin{proof}
Let ${\bf u} = u_0u_1u_2\cdots$ and ${\bf v} = v_0v_1v_2\cdots$ be any infinite
words over the alphabet $A = \{1,2,3,4\}$ that avoid squares in all arithmetic
progressions of odd difference.  We define ${\bf w}$ over the alphabet
$A \times A$ by
\[
w_{m,n} = (u_m,v_n).
\]
Consider an arbitrary line
\begin{eqnarray*}
{\bf x} & = & (w_{i_1+j_1t,i_2+j_2t})_{t \geq 0}, \\
& = & (u_{i_1+j_1t},v_{i_2+j_2t})_{t \geq 0},
\end{eqnarray*}
for some $i_1,i_2$, $j_1,j_2$, with $\gcd(j_1,j_2) = 1$.  Without loss
of generality, we may assume $j_1$ is odd.  Then the word
$(u_{i_1+j_1t})_{t \geq 0}$ is an arithmetic subsequence of odd difference
of ${\bf u}$ and hence is squarefree.  The line ${\bf x}$ is therefore
also squarefree.
\end{proof}

A backtracking search shows that there are no $2$-dimensional words ${\bf w}$
over a $7$-letter alphabet, such that every line of ${\bf w}$ is squarefree.
It remains an open problem to determine if the alphabet size of $16$ in
Theorem~\ref{carpi2} is best possible.

\begin{figure}[p]
\centering
\includegraphics{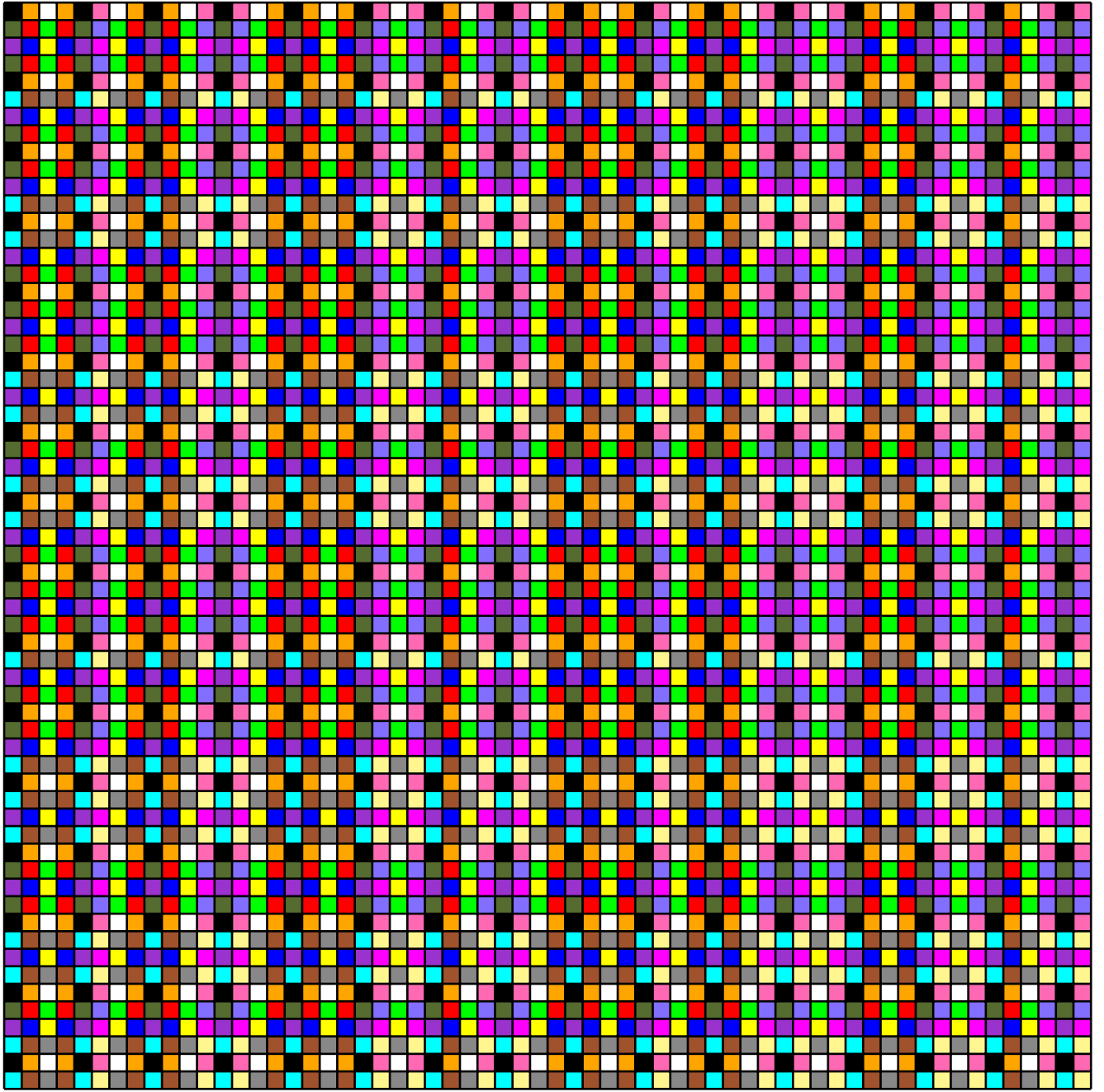}
\caption{A tiling of the $2$-dimensional grid given by a word ${\bf w}$
of Theorem~\ref{carpi2}}
\label{carpi16}
\end{figure}

Figure~\ref{carpi16} shows a tiling of the $2$-dimensional
grid induced by a word ${\bf w}$ of Theorem~\ref{carpi2}.  The colour of
the grid cell in position $(i,j)$ is determined by the value of $w_{i,j}$.

Using the results of Theorems~\ref{threeplus},~\ref{overlaps},
and~\ref{largesqs} respectively, one proves the following theorems
in a manner analogous to that of Theorem~\ref{carpi2}.

\begin{theorem}
\label{2d3plus}
There exists a 2-dimensional word ${\bf w}$ over a 4-letter alphabet,
such that every line of ${\bf w}$ is $3^+$-power-free.
\end{theorem}

\begin{theorem}
\label{2d2plus}
There exists a 2-dimensional word ${\bf w}$ over a 9-letter alphabet,
such that every line of ${\bf w}$ is $2^+$-power-free (overlapfree).
\end{theorem}

\begin{theorem}
\label{2dgeq3}
There exists a 2-dimensional word ${\bf w}$ over a 4-letter alphabet,
such that every line of ${\bf w}$ avoids squares $xx$, where $|x| \geq 3$.
\end{theorem}

\begin{figure}[p]
\centering
\includegraphics{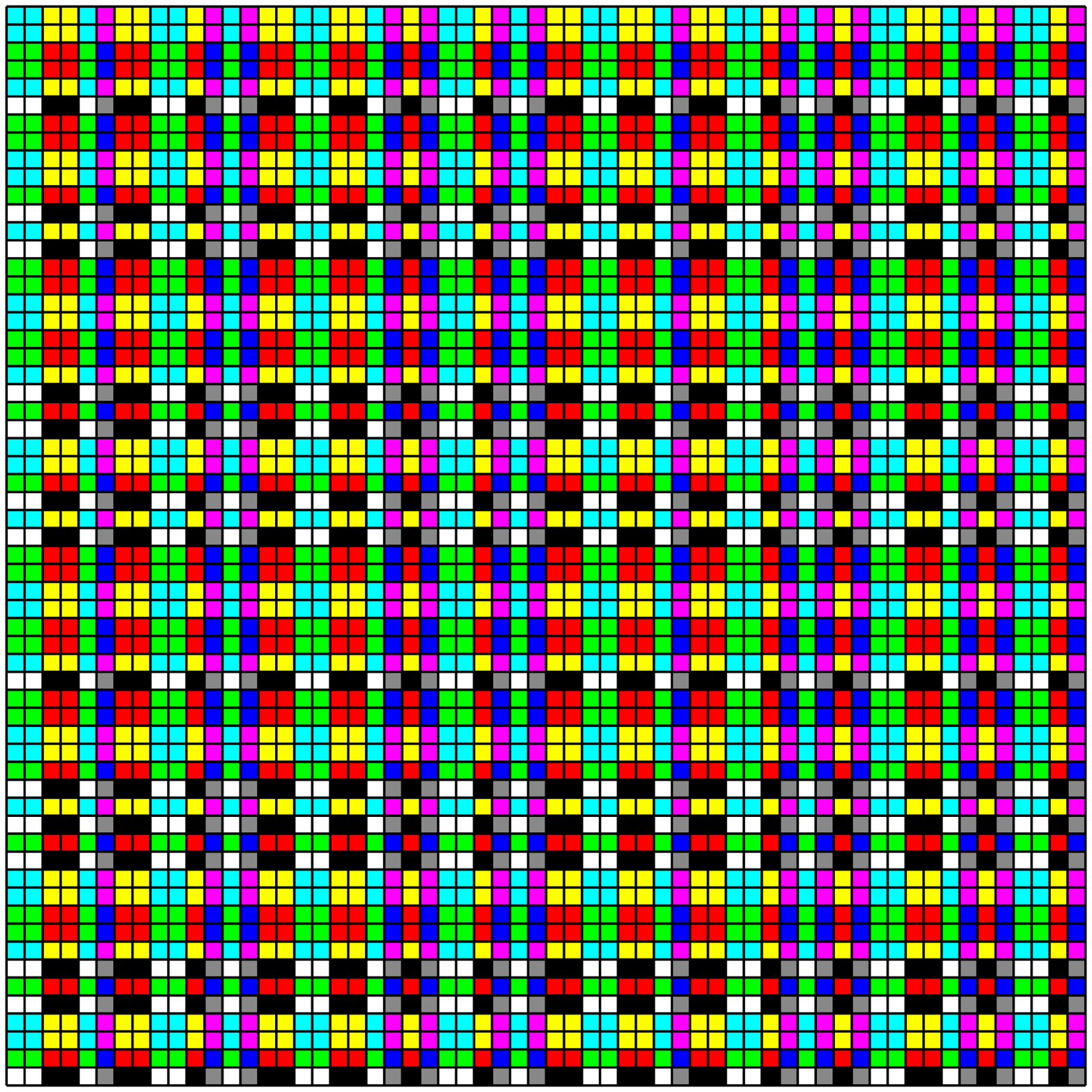}
\caption{A tiling of the $2$-dimensional grid given by a word ${\bf w}$
of Theorem~\ref{2d2plus}}
\label{carpi-overlap}
\end{figure}

\begin{figure}[p]
\centering
\includegraphics{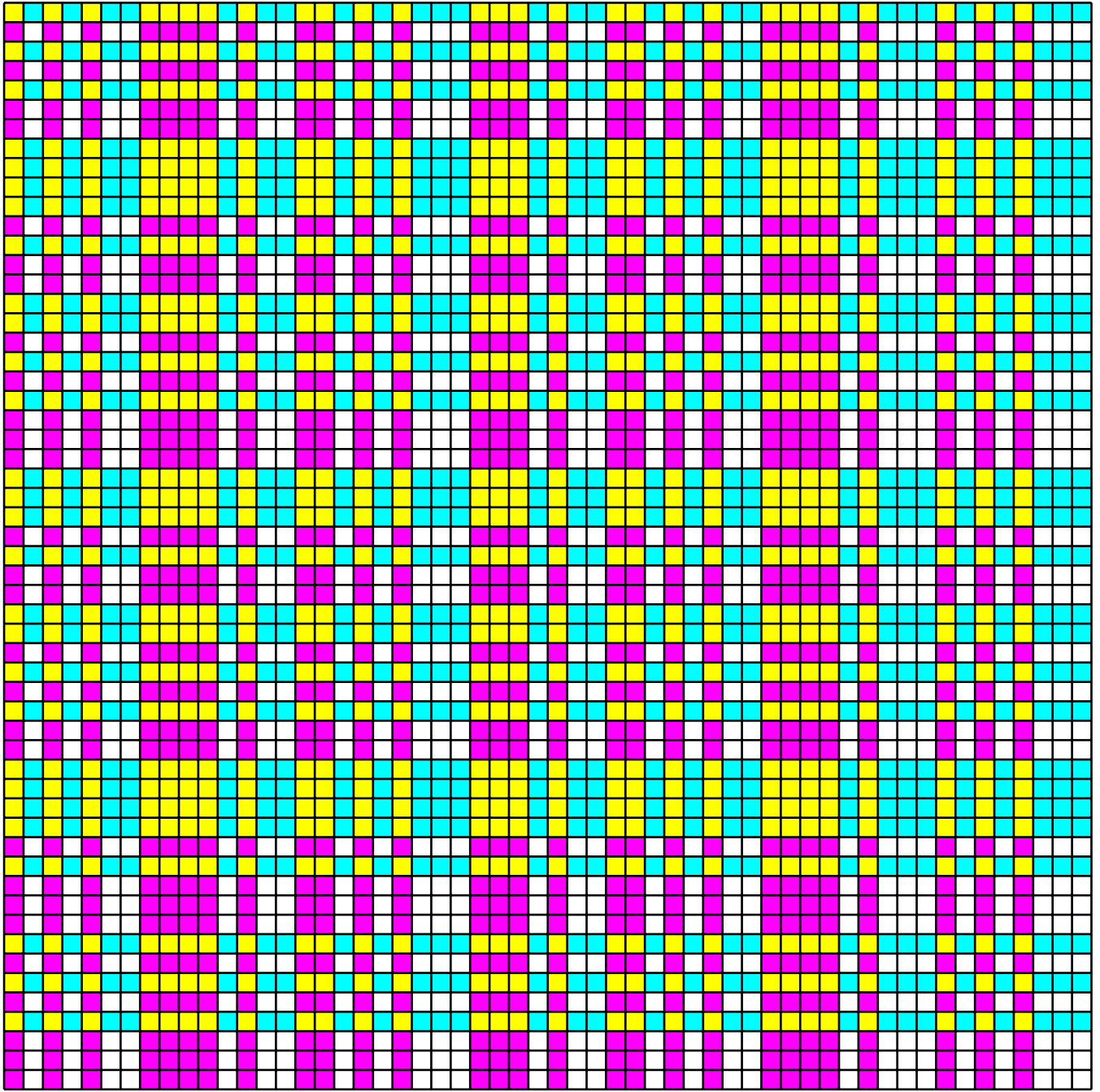}
\caption{A tiling of the $2$-dimensional grid given by a word ${\bf w}$
of Theorem~\ref{2dgeq3}}
\label{carpi-design}
\end{figure}

The reader will easily see how to generalize these results to higher
dimensions.  Figures~\ref{carpi-overlap} and \ref{carpi-design} show
tilings of the $2$-dimensional grid induced by words ${\bf w}$ of
Theorems~\ref{2d2plus} and \ref{2dgeq3}, respectively.

Grytczuk \cite{Gry06} presented the problem of determining the
\emph{Thue threshold} of $\mathbb{N}^2$, namely, the smallest
integer $t$ such that there exists an integer $k \geq 2$ and a
2-dimensional word ${\bf w}$ over a $t$-letter alphabet such that
every line of ${\bf w}$ is $k$-power-free.  Carpi's result showed
that $t \leq 16$; Theorem~\ref{2d3plus} shows that $t \leq 4$.

\section{Acknowledgments}

The second and third authors would like to thank Anna Frid for helpful
discussions.

\end{document}